\documentclass[12pt]{article}

\usepackage{amsmath}
\usepackage{amsthm}
\usepackage{amssymb}
\usepackage{latexsym}
\usepackage{pstricks}
\usepackage{graphicx, pst-plot, pst-node, pst-text, pst-tree}
\usepackage{titlefoot}
\usepackage{enumerate}

\usepackage{color}
\definecolor{lightgray}{rgb}{0.9, 0.9, 0.9}
\definecolor{darkgray}{rgb}{0.7, 0.7, 0.7}
\definecolor{darkblue}{rgb}{0, 0, .4}

\usepackage[bookmarks,breaklinks]{hyperref}
\hypersetup{
	colorlinks=true,
	linkcolor=darkblue,
	anchorcolor=darkblue,
	citecolor=darkblue,
	pagecolor=darkblue,
	urlcolor=darkblue,
	pdfpagemode=UseThumbs,
	pdftitle={Decomposing simple permutations, with enumerative consequences},
	pdfsubject={Combinatorics},
	pdfauthor={Robert Brignall, Sophie Huczynska, and Vincent Vatter},
	pdfkeywords={permutation class, restricted permutation, simple permutation}
}

\newtheorem{theorem}{Theorem}[section]
\newtheorem{proposition}[theorem]{Proposition}
\newtheorem{lemma}[theorem]{Lemma}

\newtheorem{corollary}[theorem]{Corollary}

\newcounter{todocounter}

\makeatletter
\newfont{\footsc}{cmcsc10 at 8truept}
\newfont{\footbf}{cmbx10 at 8truept}
\newfont{\footrm}{cmr10 at 10truept}
\datefoot{\today}
\amssubj{05D99, 05A15, 06A07}

\title{Decomposing simple permutations, with enumerative consequences}
\author{Robert Brignall, Sophie Huczynska\thanks{Supported by a Royal Society Dorothy Hodgkin Research Fellowship.}, and Vincent Vatter\thanks{Supported by EPSRC grant GR/S53503/01.}\\
\small School of Mathematics and Statistics\\[-5pt]
\small University of St Andrews\\[-5pt]
\small St Andrews, Fife, Scotland\\[-5pt]
\small \texttt{\{robertb, sophieh, vince\}@mcs.st-and.ac.uk}\\[-5pt]
\small \texttt{http://turnbull.mcs.st-and.ac.uk/\~{}\{\href{http://turnbull.mcs.st-and.ac.uk/~robertb}{robertb}, \href{http://turnbull.mcs.st-and.ac.uk/~sophieh}{sophieh}, \href{http://turnbull.mcs.st-and.ac.uk/~vince}{vince}\}}\\[-10pt]}

\date{}

\begin{document}
\maketitle

\newcommand{\Av}{\operatorname{Av}}
\newcommand{\C}{\mathcal{C}}
\newcommand{\W}{\mathcal{W}}
\newcommand{\Si}{\operatorname{Si}}
\newcommand{\Dash}{{\mbox{\small-}}}
\newcommand{\rect}{\operatorname{rect}}
\newcommand{\strongcomp}{\operatorname{sc}}
\newcommand{\minisec}[1]{\bigskip\noindent{\bf #1.}}

\begin{abstract}
We prove that every sufficiently long simple permutation contains two long almost disjoint simple subsequences.  This result has applications to the enumeration of restricted permutations.  For example, it immediately implies a result of B\'ona and (independently) Mansour and Vainshtein that for any $r$, the number of permutations with at most $r$ copies of $132$ has an algebraic generating function.
\end{abstract}

\section{Statement of theorem}\label{sp2-intro}

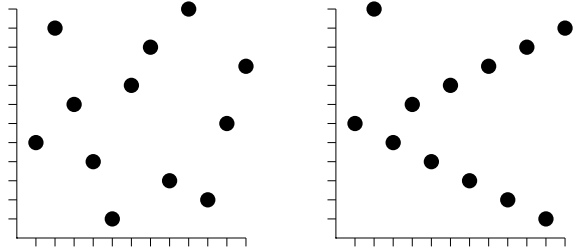
\begin{figure}
\begin{center}
\begin{tabular}{ccc}
\psset{xunit=0.01in, yunit=0.01in}
\psset{linewidth=0.005in}
\begin{pspicture}(0,0)(120,120)
\psaxes[dy=10,Dy=1,dx=10,Dx=1,tickstyle=bottom,showorigin=false,labels=none](0,0)(120,120)
\pscircle*(10,50){0.04in}
\pscircle*(20,110){0.04in}
\pscircle*(30,70){0.04in}
\pscircle*(40,40){0.04in}
\pscircle*(50,10){0.04in}
\pscircle*(60,80){0.04in}
\pscircle*(70,100){0.04in}
\pscircle*(80,30){0.04in}
\pscircle*(90,120){0.04in}
\pscircle*(100,20){0.04in}
\pscircle*(110,60){0.04in}
\pscircle*(120,90){0.04in}
\end{pspicture}
&\rule{10pt}{0pt}&
\psset{xunit=0.01in, yunit=0.01in}
\psset{linewidth=0.005in}
\begin{pspicture}(0,0)(120,120)
\psaxes[dy=10, Dy=1, dx=10, Dx=1, tickstyle=bottom, showorigin=false, labels=none](0,0)(120,120)
\pscircle*(10,60){0.04in}
\pscircle*(20,120){0.04in}
\pscircle*(30,50){0.04in}
\pscircle*(40,70){0.04in}
\pscircle*(50,40){0.04in}
\pscircle*(60,80){0.04in}
\pscircle*(70,30){0.04in}
\pscircle*(80,90){0.04in}
\pscircle*(90,20){0.04in}
\pscircle*(100,100){0.04in}
\pscircle*(110,10){0.04in}
\pscircle*(120,110){0.04in}
\end{pspicture}
\end{tabular}
\end{center}
\caption{The plots of two simple permutations.  Note that every simple subsequence of the permutation on the right must contain its first two entries.}\label{fig-wedge-plus}
\end{figure}

Simplicity, under a variety of names%
\footnote{Two synonyms for simplicity are primality and indecomposability.}, has been studied for a wide range of combinatorial objects.  Our main result concerns simple permutations; possible analogues for other contexts are discussed in the conclusion.  An {\it interval\/} in the permutation $\pi$ is a set of contiguous indices $I=[a,b]$ such that the set of values $\pi(I)=\{\pi(i) : i\in I\}$ also forms an interval of natural numbers.  Every permutation $\pi$ of $[n]=\{1,2,\dots,n\}$ has intervals of size $0$, $1$, and $n$; $\pi$ is said to be {\it simple\/} if it has no other intervals.  Figure~\ref{fig-wedge-plus} shows the plots of two simple permutations.  Intervals of permutations are interesting in their own right and have applications to biomathematics; see Corteel, Louchard, and Pemantle~\cite{corteel:common-interval:}, where among other results it is proved that the number of simple permutations of $[n]$ is asymptotic to $n!/e^2$.  More precise asymptotics are given by Albert, Atkinson, and Klazar~\cite{albert:the-enumeration:}.

Each sequence of distinct real numbers is {\it order isomorphic\/} to a unique permutation; this is the permutation with the same relative comparisons.  We say that a sequence of distinct real numbers is simple if it is order isomorphic to a simple permutation.  We prove that long simple permutations must contain two long almost disjoint simple subsequences.  Formally:

\begin{theorem}\label{sp2-main}
There is a function $f(k)$ such that every simple permutation of length at least $f(k)$ contains two simple subsequences, each of length at least $k$, sharing at most two entries.
\end{theorem}

The second ``two'' in the statement of Theorem~\ref{sp2-main} is best possible, as is demonstrated by the family of simple permutations of the form $m(2m)(m-1)(m+1)(m-2)(m+2)\cdots 1(2m-1)$; the permutation on the right of Figure~\ref{fig-wedge-plus} is of this form.  On the other hand, no attempt has been made to optimise the function $f$; our proof gives an $f$ of order about $k^{k^k}$.  The implications of Theorem~\ref{sp2-main} are discussed in the next section.  The proof begins in Section 3.

\section{Implications/motivation}\label{sp2-implications}

The permutation $\pi$ is said to {\it contain\/} the permutation $\sigma$, written $\sigma\le\pi$, if $\pi$ has a subsequence that is order isomorphic to $\sigma$.  For example, $\pi=391867452$ contains $\sigma=51342$, as can be seen by considering the subsequence $91672$ ($=\pi(2),\pi(3),\pi(5),\pi(6),\pi(9)$), and such a subsequence is called a {\it copy\/} of $\sigma$ in $\pi$.  This pattern-containment relation is a partial order on permutations.  We refer to downsets of permutations under this order as {\it permutation classes\/}.  In other words, if $\C$ is a permutation class, $\pi\in\C$, and $\sigma\le\pi$, then $\sigma\in\C$.
We denote by $\C_n$ the set $\C \cap S_n$, i.e.\ the permutations in $\C$ of length
$n$, and we refer to $\sum |\C_n| x^n$ as the {\it generating function for $\C$\/}.
Recall that
an {\it antichain} is a set of pairwise incomparable elements.
For any permutation class $\C$, there is a unique (possibly infinite) antichain $B$ such that $\C=\Av(B)=\{\pi: \beta \not \leq\pi\mbox{ for all } \beta \in B\}$. This antichain $B$, which consists of the minimal permutations not in $\C$, is called the {\it basis} of $\C$.

In a class with only finitely many simple permutations, long permutations must map nontrivial intervals onto intervals.  Thus these classes have a recursive structure in which long permutations are built up from smaller permutations, and so it is natural to expect them to have algebraic generating functions.  This is indeed the case:

\begin{theorem}[Albert and Atkinson~\cite{albert:simple-permutat:}]\label{simple}
A permutation class with only finitely many simple permutations has a readily computable algebraic generating function.
\end{theorem}

One of the simplest classes with only finitely many simple permutations is $\Av(132)$%
\footnote{In any permutation from $\Av(132)$, all entries to the left of the maximum must be greater than all entries to the right.  This shows that $\Av(132)$ has only three simple permutations ($1$, $12$, and $21$).}.  Theorems~\ref{sp2-main} and \ref{simple} combine to give a short proof of the following result.

\begin{theorem}[B\'ona~\cite{bona:the-number-of-p:}; Mansour and Vainshtein~\cite{mansour:counting-occurr:}]\label{copies-of-132}
For every $r$, the class of all permutations containing at most $r$ copies of $132$ has an algebraic generating function%
\footnote{For example, the generating function in the $r=1$ case is
$$
\frac{1-\sqrt{1-4x}}{2x}+\frac{8x^3}{\sqrt{1-4x}\left(1+\sqrt{1-4x}\right)^3}
$$
(due, originally, to B\'ona~\cite{bona:permutations-wi:}).}.
\end{theorem}

\noindent{\it Proof of Theorem~\ref{copies-of-132} via Theorems~\ref{sp2-main} and \ref{simple}.\/}
We wish to show that only finitely many simple permutations contain at most $r$ copies of $132$, or in other words, that there is a function $g(r)$ so that every simple permutation of length at least $g(r)$ contains more than $r$ copies of $132$.  Footnote 2 shows that we may take $g(0)=3$.  We now proceed by induction, setting $g(r)=f(g(\lfloor r/2\rfloor))$, where $f$ is the function from Theorem~\ref{sp2-main}.  By that theorem, every simple permutation $\pi$ of length at least $g(r)$ contains two simple subsequences of length at least $g(\lfloor r/2 \rfloor)$.  By induction each of these simple subsequences contains more than $\lfloor r/2 \rfloor $ copies of $132$.  Moreover, because these simple subsequences share at most two entries, their copies of $132$ are distinct, and thus $\pi$ contains more than $r$ copies of $132$, as desired.\qed\bigskip

Indeed, the proof above shows that every permutation class whose members contain a bounded number of copies of $132$ has an algebraic generating function, whereas Theorem~\ref{copies-of-132} is concerned only with the entire class of permutations with at most $r$ copies of $132$.  There is of course nothing special about $132$.  Denote by
$\Av(\beta_1^{\le r_1},\beta_2^{\le r_2},\dots,\beta_k^{\le r_k})$
the class of permutations that have at most $r_1$ copies of $\beta_1$, at most $r_2$ copies of $\beta_2$, and so on%
\footnote{That this is a permutation class is clear, although finding its basis may be less obvious.  An easy argument shows that the basis elements of this class have length at most $\max\{(r_i+1)|\beta_i| : i\in[k]\}$; see Atkinson~\cite{atkinson:restricted-perm:} for the details.  One such computation: $\Av(132^{\le 1})=\Av(1243$, $1342$, $1423$, $1432$, $2143$, $35142$, $354162$, $461325$, $465132)$.}.
The proof just given can be adapted to prove the following result.

\begin{corollary}\label{sp2-main-cor}
If the class $\Av(\beta_1,\beta_2,\dots,\beta_k)$ contains only finitely many simple permutations then for all choices of nonnegative integers $r_1$, $r_2$, $\dots$, and $r_k$, the class $\Av(\beta_1^{\le r_1}$, $\beta_2^{\le r_2}, \dots, \beta_k^{\le r_k})$ also contains only finitely many simple permutations.
\end{corollary}

The largest permutation class whose only simple permutations are $1$, $12$, and $21$ is the class of {\it separable permutations\/}, $\Av(2413,3142)$.  Thus as another instance of Corollary~\ref{sp2-main-cor}, we have the following.

\begin{corollary}
For all $r$ and $s$, every subclass of $\Av(2413^{\le r}, 3142^{\le s})$ contains only finitely many simple permutations and thus has an algebraic generating function.
\end{corollary}

Theorem~\ref{sp2-main} does not apply only to permutation classes.  In Brignall, Huczynska, and Vatter~\cite{brignall:simple-permutat:}, Theorem~\ref{simple} is extended to ``finite query-complete sets of properties''.  As a specialisation of that theorem, we have the following.

\begin{theorem}[Brignall, Huczynska, and Vatter~\cite{brignall:simple-permutat:}]\label{sp1-main}
In a permutation class $\C$ with only finitely many simple permutations, the following sequences have  algebraic generating functions:
\begin{itemize}
\item the number of even permutations in $\C_n$,
\item the number of involutions in $\C_n$,
\item the number of even involutions in $\C_n$,
\item the number of alternating permutations in $\C_n$,
\item the number of permutations in $\C_n$ avoiding a finite set of blocked permutations%
\footnote{{\it Blocked permutations\/}, introduced by Babson and Steingr{\'{\i}}msson~\cite{babson:generalized-per:}, are permutations containing dashes indicating the entries that need not occur consecutively.  For example, $51342$ contains two copies of $3\Dash12$: $513$ and $534$, but note that $514$ is not a copy of $3\Dash 12$ because the $1$ and $4$ are not adjacent.}.
\end{itemize}
\end{theorem}

There are several results in the literature that follow from the combination of Theorems~\ref{sp2-main}, and \ref{sp1-main}:
\begin{itemize}
\item Even permutations in $\Av(132^{\le r})$ --- Mansour~\cite{mansour:counting-occurr:a}.
\end{itemize}
(When counting even permutations, unlike when counting all permutations, symmetry considerations reduce us to three cases of length three permutations --- $123$, $132$, and $231$ --- not two\footnote{We have thus far ignored the other case; $\Av(123)$, and thus $\Av(123^{\le r})$, contains infinitely many simple permutations, so these methods do not apply.  The class $\Av(123)$ is enumerated by the Catalan numbers, $\Av(123^{\le 1})$ was counted by Noonan~\cite{noonan:the-number-of-p:}, while $\Av(123^{\le 2})$ was counted by Fulmek~\cite{fulmek:enumeration-of-:}, proving a conjecture of Noonan and Zeilberger~\cite{noonan:the-enumeration:}.  No results for larger values are known, although Fulmek conjectures formulas for $r=3$ and $r=4$ and further conjectures that $\Av(123^{\le r})$ has an algebraic generating function for all $r$.}, and thus there is another result we can state: the even permutations in $\Av(231^{\le r})$ have an algebraic generating function for all $r$, although this result seems to have escaped print.)
\begin{itemize}
\item Involutions in $\Av(231^{\le r})$ --- Mansour, Yan, and Yang~\cite{mansour:counting-occurr:b}.  In the same reference: even involutions in $\Av(231^{\le r})$.  The case of involutions in $\Av(132^{\le r})$ is due to Mansour [private communication].
\item Alternating permutations in $\Av(132^{\le r})$ --- Mansour~\cite{mansour:restricted-132-:}.
\item Permutations with at most $r$ copies of the blocked permutation $13\Dash2$ --- Claesson and Mansour~\cite{claesson:counting-occurr:}\footnote{While avoiding $132$ and avoiding $13\Dash2$ are equivalent conditions, a permutation will tend to have fewer copies of $13\Dash2$.}.
\end{itemize}

\section{Pin Sequences}\label{sec-pins}

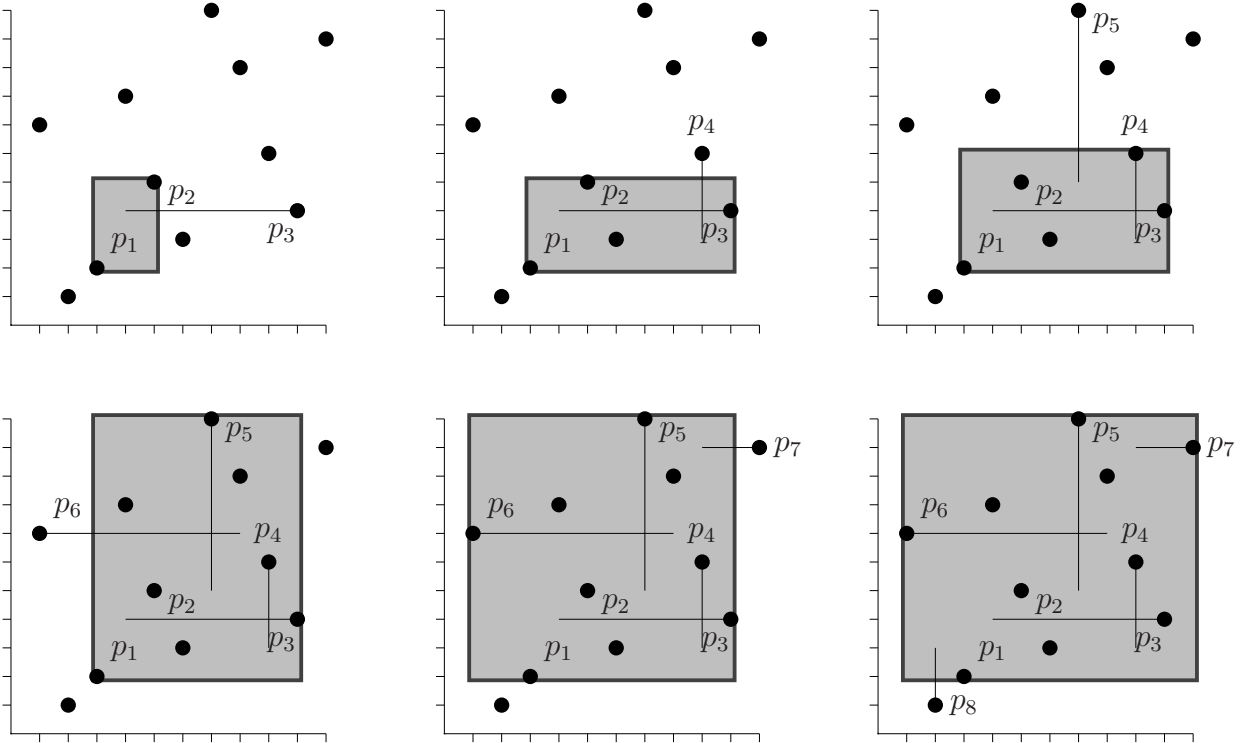
\begin{figure}
\begin{center}
\begin{tabular}{ccccc}
\psset{xunit=0.01in, yunit=0.01in}
\psset{linewidth=0.005in}
\begin{pspicture}(0,0)(180,168)
\psaxes[dy=15,dx=15,tickstyle=bottom,showorigin=false,labels=none](0,0)(165,165)
\psframe[linecolor=darkgray,fillstyle=solid,fillcolor=lightgray,linewidth=0.02in](42,27)(78,78)
\psline(150,60)(60,60)
\pscircle*(15,105){0.04in}
\pscircle*(30,15){0.04in}
\pscircle*(45,30){0.04in}
\pscircle*(60,120){0.04in}
\pscircle*(75,75){0.04in}
\pscircle*(90,45){0.04in}
\pscircle*(105,165){0.04in}
\pscircle*(120,135){0.04in}
\pscircle*(135,90){0.04in}
\pscircle*(150,60){0.04in}
\pscircle*(165,150){0.04in}
\rput[c](60,43){$p_1$}
\rput[c](90,68){$p_2$}
\rput[c](142,48){$p_3$}
\end{pspicture}
&&
\psset{xunit=0.01in, yunit=0.01in}
\psset{linewidth=0.005in}
\begin{pspicture}(0,0)(180,168)
\psaxes[dy=15,dx=15,tickstyle=bottom,showorigin=false,labels=none](0,0)(165,165)
\psframe[linecolor=darkgray,fillstyle=solid,fillcolor=lightgray,linewidth=0.02in](42,27)(153,78)
\psline(150,60)(60,60)
\psline(135,90)(135,45)
\pscircle*(15,105){0.04in}
\pscircle*(30,15){0.04in}
\pscircle*(45,30){0.04in}
\pscircle*(60,120){0.04in}
\pscircle*(75,75){0.04in}
\pscircle*(90,45){0.04in}
\pscircle*(105,165){0.04in}
\pscircle*(120,135){0.04in}
\pscircle*(135,90){0.04in}
\pscircle*(150,60){0.04in}
\pscircle*(165,150){0.04in}
\rput[c](60,43){$p_1$}
\rput[c](90,68){$p_2$}
\rput[c](142,48){$p_3$}
\rput[c](135,105){$p_4$}
\end{pspicture}
&\rule{10pt}{0pt}&
\psset{xunit=0.01in, yunit=0.01in}
\psset{linewidth=0.005in}
\begin{pspicture}(0,0)(180,168)
\psaxes[dy=15,dx=15,tickstyle=bottom,showorigin=false,labels=none](0,0)(165,165)
\psframe[linecolor=darkgray,fillstyle=solid,fillcolor=lightgray,linewidth=0.02in](42,27)(153,93)
\psline(150,60)(60,60)
\psline(135,90)(135,45)
\psline(105,165)(105,75)
\pscircle*(15,105){0.04in}
\pscircle*(30,15){0.04in}
\pscircle*(45,30){0.04in}
\pscircle*(60,120){0.04in}
\pscircle*(75,75){0.04in}
\pscircle*(90,45){0.04in}
\pscircle*(105,165){0.04in}
\pscircle*(120,135){0.04in}
\pscircle*(135,90){0.04in}
\pscircle*(150,60){0.04in}
\pscircle*(165,150){0.04in}
\rput[c](60,43){$p_1$}
\rput[c](90,68){$p_2$}
\rput[c](142,48){$p_3$}
\rput[c](135,105){$p_4$}
\rput[c](120,158){$p_5$}
\end{pspicture}
\\\\\\
\psset{xunit=0.01in, yunit=0.01in}
\psset{linewidth=0.005in}
\begin{pspicture}(0,0)(180,168)
\psaxes[dy=15,dx=15,tickstyle=bottom,showorigin=false,labels=none](0,0)(165,165)
\psframe[linecolor=darkgray,fillstyle=solid,fillcolor=lightgray,linewidth=0.02in](42,27)(153,168)
\psline(150,60)(60,60)
\psline(135,90)(135,45)
\psline(105,165)(105,75)
\psline(15,105)(120,105)
\pscircle*(15,105){0.04in}
\pscircle*(30,15){0.04in}
\pscircle*(45,30){0.04in}
\pscircle*(60,120){0.04in}
\pscircle*(75,75){0.04in}
\pscircle*(90,45){0.04in}
\pscircle*(105,165){0.04in}
\pscircle*(120,135){0.04in}
\pscircle*(135,90){0.04in}
\pscircle*(150,60){0.04in}
\pscircle*(165,150){0.04in}
\rput[c](60,43){$p_1$}
\rput[c](90,68){$p_2$}
\rput[c](142,48){$p_3$}
\rput[c](135,105){$p_4$}
\rput[c](120,158){$p_5$}
\rput[c](30,118){$p_6$}
\end{pspicture}
&\rule{10pt}{0pt}&
\psset{xunit=0.01in, yunit=0.01in}
\psset{linewidth=0.005in}
\begin{pspicture}(0,0)(180,168)
\psaxes[dy=15,dx=15,tickstyle=bottom,showorigin=false,labels=none](0,0)(165,165)
\psframe[linecolor=darkgray,fillstyle=solid,fillcolor=lightgray,linewidth=0.02in](12,27)(153,168)
\psline(150,60)(60,60)
\psline(135,90)(135,45)
\psline(105,165)(105,75)
\psline(15,105)(120,105)
\psline(165,150)(135,150)
\pscircle*(15,105){0.04in}
\pscircle*(30,15){0.04in}
\pscircle*(45,30){0.04in}
\pscircle*(60,120){0.04in}
\pscircle*(75,75){0.04in}
\pscircle*(90,45){0.04in}
\pscircle*(105,165){0.04in}
\pscircle*(120,135){0.04in}
\pscircle*(135,90){0.04in}
\pscircle*(150,60){0.04in}
\pscircle*(165,150){0.04in}
\rput[c](60,43){$p_1$}
\rput[c](90,68){$p_2$}
\rput[c](142,48){$p_3$}
\rput[c](135,105){$p_4$}
\rput[c](120,158){$p_5$}
\rput[c](30,118){$p_6$}
\rput[c](180,150){$p_7$}
\end{pspicture}
&\rule{10pt}{0pt}&
\psset{xunit=0.01in, yunit=0.01in}
\psset{linewidth=0.005in}
\begin{pspicture}(0,0)(180,168)
\psaxes[dy=15,dx=15,tickstyle=bottom,showorigin=false,labels=none](0,0)(165,165)
\psframe[linecolor=darkgray,fillstyle=solid,fillcolor=lightgray,linewidth=0.02in](12,27)(168,168)
\psline(150,60)(60,60)
\psline(135,90)(135,45)
\psline(105,165)(105,75)
\psline(15,105)(120,105)
\psline(165,150)(135,150)
\psline(30,15)(30,45)
\pscircle*(15,105){0.04in}
\pscircle*(30,15){0.04in}
\pscircle*(45,30){0.04in}
\pscircle*(60,120){0.04in}
\pscircle*(75,75){0.04in}
\pscircle*(90,45){0.04in}
\pscircle*(105,165){0.04in}
\pscircle*(120,135){0.04in}
\pscircle*(135,90){0.04in}
\pscircle*(150,60){0.04in}
\pscircle*(165,150){0.04in}
\rput[c](60,43){$p_1$}
\rput[c](90,68){$p_2$}
\rput[c](142,48){$p_3$}
\rput[c](135,105){$p_4$}
\rput[c](120,158){$p_5$}
\rput[c](30,118){$p_6$}
\rput[c](180,150){$p_7$}
\rput[c](45,15){$p_8$}
\end{pspicture}
\end{tabular}
\end{center}
\caption{A pin sequence.}
\label{fig-pins-first}
\end{figure}

Given points $p_1,\dots,p_m$ in the plane, we denote by $\rect(p_1,\dots,p_m)$ the smallest axes-parallel rectangle containing them.

Take $\pi\in S_n$ and choose two points $p_1$ and $p_2$ in the plot of $\pi$.  If these two points do not form an interval then there is at least one point which lies outside $\rect(p_1,p_2)$ and slices $\rect(p_1,p_2)$ either horizontally or vertically.  (This discussion is accompanied by the sequence of diagrams shown in Figure~\ref{fig-pins-first}.)  We call such a point a {\it pin\/}.  Choose a pin and label it $p_3$.  Now consider the larger rectangle $\rect(p_1,p_2,p_3)$.  If this also does not form an interval in $\pi$ then we can find another pin, $p_4$, which slices $\rect(p_1,p_2,p_3)$ either horizontally or vertically.  Again, if $\rect(p_1,p_2,p_3,p_4)$ is not an interval then we can find another pin $p_5$.  We refer to a sequence of pins constructed in this manner as a {\it pin sequence\/}.  

Formally, a {\it pin sequence\/} is a sequence of points $p_1$, $p_2$, $\dots$ in the plot of $\pi$ such that for each $i\ge 3$,
\begin{itemize}
\item $p_i\not\in\rect(p_1,\dots,p_{i-1})$, and
\item if $\rect(p_1, \dots, p_{i-1})=[a,b]\times[c,d]$ and $p_i=(x,y)$, we have either $a<x<b$ or $c<y<d$, or, in other words, $p_i$ {\it slices\/} $\rect(p_1,\dots,p_{i-1})$ either horizontally or vertically.
\end{itemize}
We describe pins as either {\it left\/}, {\it right\/}, {\it up\/}, or {\it down\/} based on their position relative to the rectangle that they slice.  Thus in the pin sequence from Figure~\ref{fig-pins-first}, $p_3$ and $p_7$ are right pins, $p_4$ and $p_5$ are up pins, $p_6$ is a left pin, and $p_8$ is a down pin ($p_1$ and $p_2$ lack direction).

A {\it proper pin sequence\/} is one that satisfies two additional conditions:
\begin{itemize}
\item {\it Maximality condition\/}: each pin must be maximal in its direction.  For example, if $\rect(p_1,\dots,p_{i-1})=[a,b]\times[c,d]$ and $p_i=(x,y)$ is a right pin, then it is the right-most of all possible right pins for this rectangle, or, in other words, the region $(x,n]\times [c,d]$ is devoid of points.
\item {\it Separation condition\/}: $p_{i+1}$ must {\it separate\/} $p_i$ from $\{p_1,\dots,p_{i-1}\}$.  That is, $p_{i+1}$ must lie horizontally or vertically between $\rect(p_1,\dots,p_{i-1})$ and $p_i$.
\end{itemize}
For example, in the pin sequence shown in Figure~\ref{fig-pins-first}, the choice of $p_4$ violates the maximality condition, while the choices of $p_5$, $p_7$, and $p_8$ violate the separation condition.  The ultimate goal of the following succession of lemmas is to show (in Theorem~\ref{pins-simple}) that all or all but one of the pins in a proper pin sequence themselves form a simple permutation.  We begin by observing that proper pin sequences travel by $90^{\circ}$ turns only.

\begin{lemma}\label{pinseq1}
In a proper pin sequence, $p_{i+1}$ cannot lie in the same or opposite direction as $p_i$ (for all $i\ge 3$).
\end{lemma}
\begin{proof}
By the maximality condition, $p_{i+1}$ cannot lie in  the same direction as $p_i$.  It cannot lie in the opposite direction by the separation condition.
\end{proof}

\begin{lemma}\label{pinseq2}
In a proper pin sequence, $p_i$ does not separate any two members of $\{p_1,\dots,p_{i-2}\}$.
\end{lemma}
\begin{proof}
If $p_i$ did separate $\rect(p_1,\dots,p_{i-2})$ into two parts then $p_{i-1}$ would lie on one side of this divide, violating the separation condition.
\end{proof}

\begin{lemma}\label{pinseq3}
In a proper pin sequence, $p_i$ and $p_{i+1}$ are separated either by $p_{i-1}$ or by each of $p_1, \dots, p_{i-2}$.
\end{lemma}
\begin{proof}
The lemma is vacuously true for $i=1$ and $i=2$, so let us assume that $i\ge 3$.  Without loss we may assume that $p_{i-1}$ is a right pin and $p_i$ is an up pin.  By Lemma~\ref{pinseq1}, $p_{i+1}$ must be either a right pin or a left pin.  The remainder of the proof is evident from Figure~\ref{fig-pinseq3}.
\begin{figure}
\begin{center}
\begin{tabular}{ccc}
\psset{xunit=0.015in, yunit=0.015in}
\psset{linewidth=0.005in}
\begin{pspicture}(0,-10)(145,90)
\psaxes[dy=1000, Dy=1, dx=1000, Dx=1, tickstyle=bottom, showorigin=false, labels=none](0,-10)(140,70)
\psframe[linecolor=darkgray,fillstyle=solid,fillcolor=lightgray,linewidth=0.02in](20,10)(120,50)
\rput[c](70,20){$\rect(p_1,\dots,p_{i-1})$}
\psline(110,40)(140,40)
\pscircle*(140,40){0.04in}
\rput[l](145,40){$p_{i-1}$}
\psline(130,30)(130,70)
\pscircle*(130,70){0.04in}
\rput[c](130,78){$p_{i}$}
\psline(140,60)(10,60)
\pscircle*(10,60){0.04in}
\rput[c](15,68){$p_{i+1}$}
\end{pspicture}
&\rule{10pt}{0pt}&
\psset{xunit=0.015in, yunit=0.015in}
\psset{linewidth=0.005in}
\begin{pspicture}(0,-10)(165,90)
\psaxes[dy=1000, Dy=1, dx=1000, Dx=1, tickstyle=bottom, showorigin=false, labels=none](0,-10)(160,70)
\psframe[linecolor=darkgray,fillstyle=solid,fillcolor=lightgray,linewidth=0.02in](20,10)(120,50)
\rput[c](70,20){$\rect(p_1,\dots,p_{i-1})$}
\psline(110,40)(140,40)
\pscircle*(140,40){0.04in}
\rput[l](145,40){$p_{i-1}$}
\psline(130,30)(130,70)
\pscircle*(130,70){0.04in}
\rput[c](130,78){$p_{i}$}
\psline(120,60)(160,60)
\pscircle*(160,60){0.04in}
\rput[c](160,68){$p_{i+1}$}
\end{pspicture}
\end{tabular}
\end{center}
\caption{The two cases in the proof of Lemma~\ref{pinseq3}.}
\label{fig-pinseq3}
\end{figure}
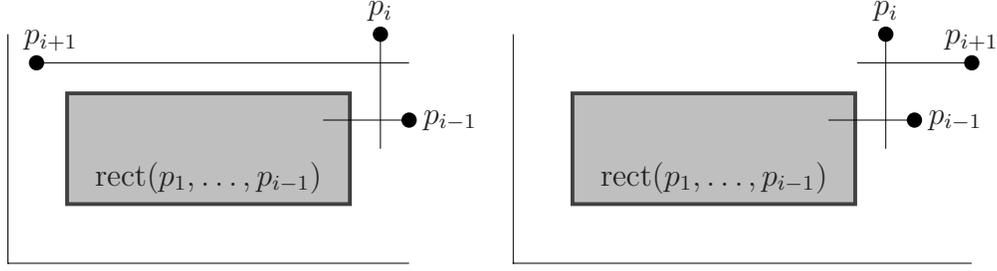
\end{proof}

We are now ready to prove our main result about proper pin sequences.

\begin{theorem}\label{pins-simple}
If $p_1,\dots,p_m$ is a proper pin sequence then one of the sets of points $\{p_1,\dots,p_m\}$, $\{p_1,\dots,p_m\}\setminus\{p_1\}$, or $\{p_1,\dots,p_m\}\setminus\{p_2\}$ is order isomorphic to a simple permutation.
\end{theorem}
\begin{proof}
Suppose $m\ge 4$, as the smaller cases are trivially true.  We are interested in the possible intervals in the subsequence given by the pins $p_1,\dots,p_m$; we shall call these {\it intervals of pins\/}.  Take $M\subseteq\{p_1,\dots,p_m\}$ to be a minimal non-singleton interval of pins.  Note that $M$ is therefore order isomorphic to a simple permutation.

If $M$ contains a pair of pins $p_i$ and $p_j$ with $i<j<m$ then by the separation condition $p_{j+1},\dots,p_m\in M$.  Furthermore, because $j<m$, Lemma~\ref{pinseq3} shows that $M$ contains either $p_{j-1}$ or $p_1,p_2,\dots,p_{j-2}$.  In the latter case, if $j\ge 4$ then separation gives $p_{j-1}\in M$, as desired, while if $j\le 3$, we have already found a minimal non-singleton interval of pins of the desired form.  In the former case, the proof is completed by iterating this process.

Only the case $M=\{p_i,p_m\}$ remains.  If $m-1$, then Lemma~\ref{pinseq3} gives a contradiction.  If $3\le i\le m-2$ then, by the separation condition, $p_i$ separates $\{p_1,\dots,p_{i-1}\}$, while Lemma~\ref{pinseq3} shows that $p_m$ does not separate these points; thus at least one of them must lie in $M$, another contradiction.

We are now reduced to the cases $M=\{p_1,p_m\}$ and $M=\{p_2,p_m\}$.  We consider the former; the latter is analogous.  Because $p_3$ separates $p_2$ from $p_1$, it also separates $p_2$ from $p_m$, so $\{p_2,p_m\}$ cannot be an interval.  If there are any other minimal non-singleton intervals of pins, then we are done by the considerations above.  Therefore, $\{p_1,p_m\}$ is the only minimal non-singleton interval of pins, and thus $\{p_2,\dots,p_m\}$ is order isomorphic to a simple permutation.
\end{proof}

As a corollary of this theorem, we see that Theorem~\ref{sp2-main} (in fact, a stronger result) is true for simple permutations with long pin sequences.

\begin{corollary}\label{long-pins}
If $\pi$ contains a proper pin sequence of length at least $2k+2$ then $\pi$ contains two disjoint simple subsequences, each of length at least $k$.
\end{corollary}
\begin{proof}
Apply Theorem~\ref{pins-simple} to the two pin sequences $p_1,\dots,p_{k+1}$ and $p_{k+2},\dots,p_{2k+2}$.
\end{proof}

We say that the pin sequence $p_1,\dots,p_m$ for the permutation $\pi\in S_n$ is {\it saturated\/} if $\rect(p_1,\dots,p_m)=[n]\times[n]$.  For example, the pin sequence in Figure~\ref{fig-pins-first} is saturated.  Any two points $p_1\neq p_2$ in the plot of a simple permutation can be extended to a saturated pin sequence, as we are forced to stop extending a pin sequence only upon finding an interval or when the rectangle contains every point in $\pi$.

It is important to note that two points in a simple permutations need not be extendable to a proper saturated pin sequence.  For example, the permutation in Figure~\ref{fig-pins-first} does not have a proper saturated pin sequence beginning with $p_1$ and $p_2$.  For this reason we work with a weaker requirement: the pin sequence $p_1,\dots,p_m$ is said to be {\it right-reaching\/} if $p_m$ is the right-most point of $\pi$.

\begin{lemma}\label{right-reaching}
For every simple permutation $\pi$ and pair of points $p_1$ and $p_2$ (unless, trivially, $p_1$ is the right-most point of $\pi$), there is a proper right-reaching pin sequence beginning with $p_1$ and $p_2$.
\end{lemma}
\begin{proof}
Clearly we can find a saturated pin sequence $p_1,p_2,\dots$ in $\pi$ that satisfies the maximality condition.  Since this pin sequence is saturated, it includes the right-most point; label it $p_{i_1}$.  Now take $i_2$ as small as possible so that $p_1, p_2, \dots, p_{i_2}, p_{i_1}$ is a valid pin sequence.  Note first that $i_2<i_1$ because $p_1,\dots,p_{i_1}$ is a valid pin sequence.  Now observe that $p_{i_1}$ separates $p_{i_2}$ from $\rect(p_1,\dots,p_{i_2-1})$, because $p_1,\dots,p_{i_2-1},p_{i_1}$ is not a valid pin sequence.  Continuing in this manner, we find pins $p_{i_3}$, $p_{i_4}$, and so on, until we reach the stage where $p_{i_{m+1}}=p_2$.  Then $p_1,p_2,p_{i_m},p_{i_{m-1}},\dots,p_{i_1}$ is a proper right-reaching pin sequence.
\end{proof}

\section{Simple permutations without long proper pin sequences}

\begin{figure}
\begin{center}
\begin{tabular}{ccc}
\psset{xunit=0.01in, yunit=0.01in}
\psset{linewidth=0.005in}
\begin{pspicture}(0,0)(150,150)
\psaxes[dy=10,Dy=1,dx=10,Dx=1,tickstyle=bottom,showorigin=false,labels=none](0,0)(140,140)
\psline[linecolor=darkgray,linewidth=0.02in](0,75)(150,75)
\pscircle*(10,70){0.04in}
\pscircle*(20,100){0.04in} \pscircle*(30,30){0.04in}
\pscircle*(40,120){0.04in} \pscircle*(50,50){0.04in}
\pscircle*(60,80){0.04in} \pscircle*(70,60){0.04in}
\pscircle*(80,140){0.04in} \pscircle*(90,10){0.04in}
\pscircle*(100,110){0.04in} \pscircle*(110,40){0.04in}
\pscircle*(120,130){0.04in} \pscircle*(130,20){0.04in}
\pscircle*(140,90){0.04in}
\end{pspicture}
&\rule{10pt}{0pt}&
\psset{xunit=0.01in, yunit=0.01in}
\psset{linewidth=0.005in}
\begin{pspicture}(0,0)(150,150)
\psline[linecolor=darkgray,linewidth=0.02in](75,0)(75,150)
\psaxes[dy=10,Dy=1,dx=10,Dx=1,tickstyle=bottom,showorigin=false,labels=none](0,0)(140,140)
\pscircle*(10,90){0.04in}
\pscircle*(20,130){0.04in} \pscircle*(30,30){0.04in}
\pscircle*(40,110){0.04in} \pscircle*(50,50){0.04in}
\pscircle*(60,70){0.04in} \pscircle*(70,10){0.04in}
\pscircle*(80,60){0.04in} \pscircle*(90,140){0.04in}
\pscircle*(100,20){0.04in} \pscircle*(110,100){0.04in}
\pscircle*(120,40){0.04in} \pscircle*(130,120){0.04in}
\pscircle*(140,80){0.04in}
\end{pspicture}
\end{tabular}
\end{center}
\caption{A vertical alternation (left) and its inverse, a
horizontal alternation (right).}\label{fig-alternations}
\end{figure}
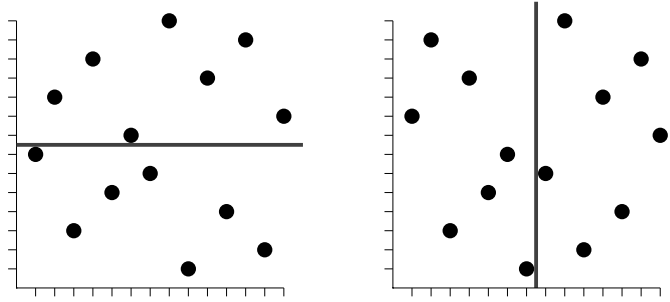

It remains only to consider simple permutations without long proper pin sequences, a consideration which constitutes the bulk of the proof.  Our goal in this section is to prove that these permutations contain long ``alternations''.  A {\it horizontal alternation\/} is a permutation in which every odd entry lies to the left of every even entry, or the reverse of such a permutation.  A {\it vertical alternation\/} is the group-theoretic inverse of a horizontal alternation.  Examples are shown in Figure~\ref{fig-alternations}.

Every sufficiently long vertical alternation contains either a long {\it parallel alternation\/} or a long {\it wedge alternation\/} (see Figure~\ref{fig-wedge-parallel} for definitions):

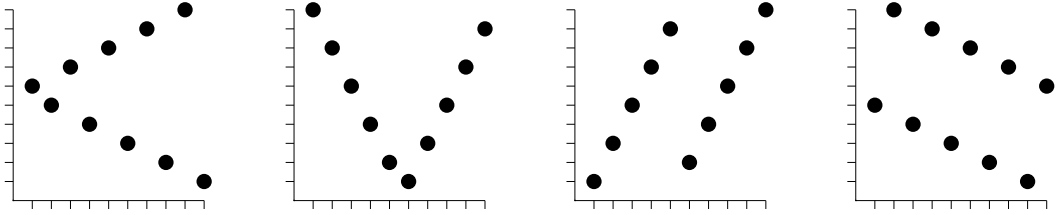
\begin{figure}
\begin{center}
\begin{tabular}{ccccccc}
\psset{xunit=0.01in, yunit=0.01in}
\psset{linewidth=0.005in}
\begin{pspicture}(0,0)(100,100)
\psaxes[dy=10, Dy=1, dx=10, Dx=1, tickstyle=bottom, showorigin=false, labels=none](0,0)(100,100)
\pscircle*(10,60){0.04in}
\pscircle*(20,50){0.04in}
\pscircle*(30,70){0.04in}
\pscircle*(40,40){0.04in}
\pscircle*(50,80){0.04in}
\pscircle*(60,30){0.04in}
\pscircle*(70,90){0.04in}
\pscircle*(80,20){0.04in}
\pscircle*(90,100){0.04in}
\pscircle*(100,10){0.04in}
\end{pspicture}
&\rule{10pt}{0pt}&
\psset{xunit=0.01in, yunit=0.01in}
\psset{linewidth=0.005in}
\begin{pspicture}(0,0)(100,100)
\psaxes[dy=10,Dy=1,dx=10,Dx=1,tickstyle=bottom,showorigin=false,labels=none](0,0)(100,100)
\pscircle*(10,100){0.04in}
\pscircle*(20,80){0.04in}
\pscircle*(30,60){0.04in}
\pscircle*(40,40){0.04in}
\pscircle*(50,20){0.04in}
\pscircle*(60,10){0.04in}
\pscircle*(70,30){0.04in}
\pscircle*(80,50){0.04in}
\pscircle*(90,70){0.04in}
\pscircle*(100,90){0.04in}
\end{pspicture}
&\rule{10pt}{0pt}&
\psset{xunit=0.01in, yunit=0.01in}
\psset{linewidth=0.005in}
\begin{pspicture}(0,0)(100,100)
\psaxes[dy=10,Dy=1,dx=10,Dx=1,tickstyle=bottom,showorigin=false,labels=none](0,0)(100,100)
\pscircle*(10,10){0.04in}
\pscircle*(20,30){0.04in}
\pscircle*(30,50){0.04in}
\pscircle*(40,70){0.04in}
\pscircle*(50,90){0.04in}
\pscircle*(60,20){0.04in}
\pscircle*(70,40){0.04in}
\pscircle*(80,60){0.04in}
\pscircle*(90,80){0.04in}
\pscircle*(100,100){0.04in}
\end{pspicture}
&\rule{10pt}{0pt}&
\psset{xunit=0.01in, yunit=0.01in}
\psset{linewidth=0.005in}
\begin{pspicture}(0,0)(100,100)
\psaxes[dy=10, Dy=1, dx=10, Dx=1, tickstyle=bottom, showorigin=false, labels=none](0,0)(100,100)
\pscircle*(10,50){0.04in}
\pscircle*(20,100){0.04in}
\pscircle*(30,40){0.04in}
\pscircle*(40,90){0.04in}
\pscircle*(50,30){0.04in}
\pscircle*(60,80){0.04in}
\pscircle*(70,20){0.04in}
\pscircle*(80,70){0.04in}
\pscircle*(90,10){0.04in}
\pscircle*(100,60){0.04in}
\end{pspicture}
\end{tabular}
\end{center}
\caption{The two permutations on the left are wedge alternations, the two on the right are parallel alternations.}
\label{fig-wedge-parallel}
\end{figure}

\begin{proposition}\label{long-alternations}
Every alternation of length at least $2k^4$ contains either a parallel or wedge alternation of length at least $2k$.
\end{proposition}
\begin{proof}
Let $\pi$ be a vertical alternation of length $2n\ge 2k^4$.  By the Erd\H os-Szekeres Theorem (every permutation of length $n$ contains a monotone subsequence of length at least $\sqrt{n}$), the sequence $\pi(1), \pi(3), \dots, \pi(2n-1)$ contains a monotone subsequence of length at least $k^2$, say $\pi(i_1),\pi(i_2),\dots,\pi(i_{k^2})$.  Applying the Erd\H os-Szekeres Theorem to the subsequence $\pi(i_1+1),\pi(i_2+1),\dots,\pi(i_{k^2}+1)$ completes the proof.
\end{proof}

Note that every parallel alternation of length $2k+2\ge 10$ contains two disjoint simple permutations of length at least $k$.  Thus Theorem~\ref{sp2-main} follows in the case where our simple permutation contains a long parallel alternation.

We say that the pin sequences $p_1,\dots,p_s$ and $q_1,\dots,q_t$ {\it converge at the pin $x$\/} if there exist $i, j\ge 3$ so that $p_i=q_j=x$ but $p_{i-1}\neq q_{j-1}$.

\begin{lemma}\label{convergence-lemma}
If $8k$ proper pin sequences of $\pi$ converge at the same pin, then $\pi$ contains an alternation of length at least $2k$.
\end{lemma}
\begin{proof}
Let us suppose that $8k$ pin sequences converge at the pin $p$.  This pin could be variously functioning as a left, right, down, or up pin for each of these $8k$ sequences, but $p$ plays the same role for at least $2k$ sequences.  Suppose, without loss, that $p$ is a right pin for at least $2k$ sequences.  Now consider the immediate predecessors to $p$ in these sequences.  These pins can be either up pins or down pins (by Lemma~\ref{pinseq1}).  By symmetry, we may assume that for at least $k$ of these pin sequences the immediate predecessor to $p$ is an up pin.  Reading left to right, label these immediate predecessor pins $p^{(1)}$, $p^{(2)}, \dots, p^{(k)}$ and let $R^{(i)}$ denote the rectangle for which $p^{(i)}$ is a pin.  Note that each $R^{(i)}$ lies completely below $p$, as otherwise the separation condition would prevent $p$ from following $p^{(i)}$ in the corresponding pin sequence.  We now have the situation depicted in Figure~\ref{fig-convergence}.

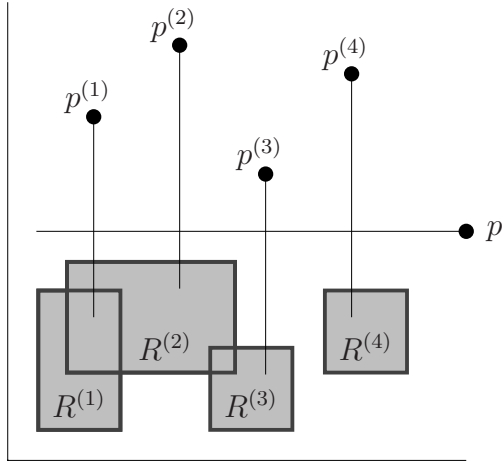
\begin{figure}
\begin{center}
\psset{xunit=0.015in, yunit=0.015in}
\psset{linewidth=0.005in}
\begin{pspicture}(0,0)(165,165)
\psaxes[dy=1000, Dy=1, dx=1000, Dx=1, tickstyle=bottom, showorigin=false, labels=none](0,0)(160,160)
\psframe[linecolor=darkgray,fillstyle=solid,fillcolor=lightgray,linewidth=0.02in](10,10)(40,60)
\psframe[linecolor=darkgray,fillstyle=solid,fillcolor=lightgray,linewidth=0.02in](20,30)(80,70)
\psframe[linecolor=darkgray,fillstyle=solid,fillcolor=lightgray,linewidth=0.02in](70,10)(100,40)
\psframe[linecolor=darkgray,fillstyle=solid,fillcolor=lightgray,linewidth=0.02in](110,30)(140,60)
\psframe[linecolor=darkgray,linewidth=0.02in](10,10)(40,60)
\psframe[linecolor=darkgray,linewidth=0.02in](20,30)(80,70)
\rput[c](25,20){$R^{(1)}$}
\rput[c](55,40){$R^{(2)}$}
\rput[c](85,20){$R^{(3)}$}
\rput[c](125,40){$R^{(4)}$}
\psline(10,80)(160,80)
\pscircle*(160,80){0.04in}
\rput[c](170,80){$p$}
\psline(30,50)(30,120)
\pscircle*(30,120){0.04in}
\rput[c](28,127){$p^{(1)}$}
\psline(60,60)(60,145)
\pscircle*(60,145){0.04in}
\rput[c](58,152){$p^{(2)}$}
\psline(90,30)(90,100)
\pscircle*(90,100){0.04in}
\rput[c](88,107){$p^{(3)}$}
\psline(120,50)(120,135)
\pscircle*(120,135){0.04in}
\rput[c](118,142){$p^{(4)}$}
\end{pspicture}
\end{center}
\caption{The situation that arises in the proof of Lemma~\ref{convergence-lemma}.}
\label{fig-convergence}
\end{figure}

It suffices to show, for each $i$, that $\pi$ contains a point lying horizontally between $p^{(i)}$ and $p^{(i+1)}$ and below $p$, since then these points, together with the $p^{(i)}$'s and $p$, will give an alternation of length $2k$.  However, if there is no such point then $p^{(i)}$ and $p^{(i+1)}$ could each function as up pins for both $R^{(i)}$ and $R^{(i+1)}$, and thus one of these choices would contradict the maximality condition, completing the proof.
\end{proof}

\begin{lemma}\label{pin-or-alt}
Every simple permutation of length at least $2(8k^4)^{(8k^4)^{2k}}$ contains either a proper pin sequence of length at least $2k$ or a parallel or wedge alternation of length at least $2k$.
\end{lemma}
\noindent{\it Proof.}
Suppose that the simple permutation $\pi\in S_n$ contains neither a proper pin sequence of length at least $2k$ nor a parallel or wedge alternation of length at least $2k$.  In particular, $\pi$ does not contain a proper right-reaching pin sequence of length $2k$, and it follows from Proposition~\ref{long-alternations} that $\pi$ has no alternations of length $2k^4$.

Pair up each of the entries of $\pi$ except the right-most.  Taking proper right-reaching pin sequences beginning at each of these pairs creates $\lfloor (n-1)/2\rfloor$ sequences.

As these pin sequences are right-reaching, they all agree on their final (right-most) pin which we denote by $p$.  By Lemma~\ref{convergence-lemma}, fewer than $8k^4$ of these pin sequences converge at $p$; equivalently, there are fewer than $8k^4$ immediate predecessors to $p$.  Label these immediate predecessors $p^{(1)},p^{(2)},\dots,p^{(m)}$.  Again, fewer than $8k^4$ pin sequences converge at each of the $p^{(i)}$'s, so there are fewer than $(8k^4)^{8k^4}$ immediate predecessors to these pins.  Continue this process until we reach the sequences of length $2k$, of which we have assumed there are none.  We have thus counted all $\lfloor (n-1)/2\rfloor$ of our sequences, and have obtained the bound
$$
\left\lfloor\frac{n-1}{2}\right\rfloor
<
1+8k^4+(8k^4)^{8k^4}+(8k^4)^{(8k^4)^2}+\dots+(8k^4)^{(8k^4)^{(2k-1)}},
$$
so, simplifying,\\
\newcommand{\eqed}[1]{$\textcolor{white}{\qed}\hfill{#1}\hfill\qed$}
\eqed{n<2(8k^4)^{(8k^4)^{2k}}.}

\begin{figure}
\begin{center}
\begin{tabular}{ccc}
\psset{xunit=0.01in, yunit=0.01in}
\psset{linewidth=0.005in}
\begin{pspicture}(0,0)(120,120)
\psaxes[dy=10, Dy=1, dx=10, Dx=1, tickstyle=bottom, showorigin=false,labels=none](0,0)(120,120)
\pscircle*(10,70){0.04in}
\pscircle*(20,50){0.04in}
\pscircle*(30,80){0.04in}
\pscircle*(40,40){0.04in}
\pscircle*(50,90){0.04in}
\pscircle*(60,30){0.04in}
\pscircle*(70,100){0.04in}
\pscircle*(80,20){0.04in}
\pscircle*(90,110){0.04in}
\pscircle*(100,10){0.04in}
\pscircle*(110,120){0.04in}
\pscircle*(120,60){0.04in}
\end{pspicture}
&\rule{10pt}{0pt}&
\psset{xunit=0.01in, yunit=0.01in}
\psset{linewidth=0.005in}
\begin{pspicture}(0,0)(120,120)
\psaxes[dy=10, Dy=1, dx=10, Dx=1, tickstyle=bottom, showorigin=false, labels=none](0,0)(120,120)
\pscircle*(10,60){0.04in}
\pscircle*(20,120){0.04in}
\pscircle*(30,50){0.04in}
\pscircle*(40,70){0.04in}
\pscircle*(50,40){0.04in}
\pscircle*(60,80){0.04in}
\pscircle*(70,30){0.04in}
\pscircle*(80,90){0.04in}
\pscircle*(90,20){0.04in}
\pscircle*(100,100){0.04in}
\pscircle*(110,10){0.04in}
\pscircle*(120,110){0.04in}
\end{pspicture}
\end{tabular}
\end{center}
\caption{The two types of wedge simple permutations.  On the left, type $1$, on the right, type $2$.}\label{fig-fund-wedge}
\end{figure}
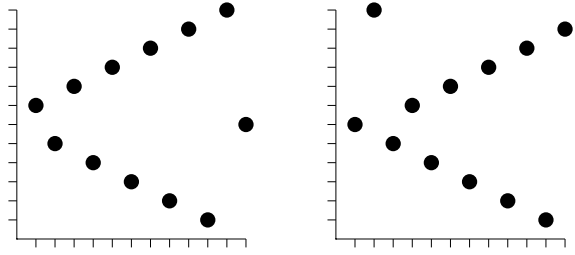

We are left to deal with simple permutations which do not have long proper pin sequences but do have long wedge alternations.  We prove that these permutations contain long {\it wedge simple permutations\/}, of which there are two types (up to symmetry).  Examples of these two types are shown in Figure~\ref{fig-fund-wedge}.

\begin{lemma}\label{long-wedge}
If a simple permutation contains a wedge alternation of length $4k^2$ then it contains either a pin sequence of length at least $2k$ or a wedge simple permutation of length at least $2k$.
\end{lemma}
\begin{proof}
Let $\pi$ be a simple permutation containing a wedge alternation of length at least $4k^2$.  By symmetry we may assume that this wedge alternation opens to the right (i.e., it is oriented as $<$).  We call these the {\it wedge points\/} of $\pi$.  Label the two left-most wedge points $p_1$ and $p_2$ and by Lemma~\ref{right-reaching} extend this into a proper right-reaching pin sequence $p_1,p_2,\dots,p_m$.

\newcommand{\ws}{\operatorname{ws}}
\newcommand{\wc}{\operatorname{wc}}

Let $R_i$ denote the smallest rectangle in the plot of $\pi$ containing $p_1$, $p_2$, and $p_i$ that is not sliced by a wedge point outside the rectangle.  Define the {\it wedge sum\/} of the pin $p_i$, $\ws(p_i)$, to be the number of wedge points in $R_i$.  For $i\ge 2$ define the {\it wedge contribution\/} of $p_i$ by $\wc(p_i)=\ws(p_i)-\ws(p_{i-1})$ and set $\wc(p_1)=1$.  Regarding these quantities we make four observations:
\begin{enumerate}[(W1)]
\item the wedge sum of $p_m$ is equal to the total number of wedge points and also to $\displaystyle\sum_{i=1}^m \wc(p_i)$,
\item it is not hard to construct examples in which pins have negative wedge contributions; indeed,
\item left pins cannot have positive wedge contributions, and finally,
\item if $p_i$ is an up pin, then the right-most wedge point in $R_i$ is an upper wedge point.
\end{enumerate}

We now claim that each $p_i$ lies in a wedge simple permutation of length at least $\wc(p_i)+2$.  This claim implies the theorem, because if no pin lies in a wedge simple permutation of length at least $2k$ then $\wc(p_i)\le 2k-3$, so by (W1),
$$
4k^2\le\sum_{i=1}^m \wc(p_i)\le m(2k-3),
$$
and thus $m\ge 2k$, giving the long pin sequence desired.

The claim is easily observed for $i=1$ and, by (W3), vacuously true if $p_i$ is a left pin.  Thus by symmetry there are only three cases to consider: an up pin followed by a right pin, a right pin followed by an up pin, and a left pin followed by an up pin.  These three cases are depicted in Figure~\ref{fig-long-wedge}.

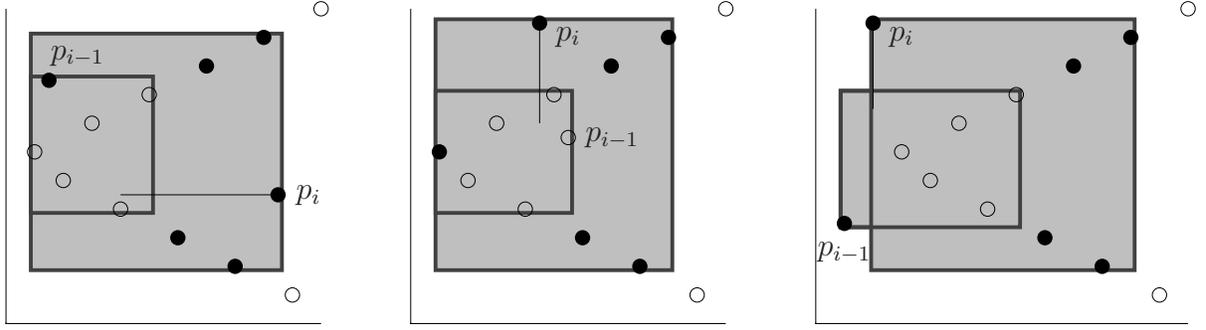
\begin{figure}
\begin{center}
\begin{tabular}{ccccc}
\psset{xunit=0.01in, yunit=0.01in}
\psset{linewidth=0.005in}
\begin{pspicture}(0,0)(165,165)
\psaxes[dy=1000,dx=1000](0,0)(165,165)
\psframe[linecolor=darkgray,fillstyle=solid,fillcolor=lightgray,linewidth=0.02in](12,27)(145.5,153)
\psframe[linecolor=darkgray,fillstyle=solid,fillcolor=lightgray,linewidth=0.02in](12,57)(78,130.5)
\psline(142.5,67.5)(60,67.5)
\pscircle*(22.5,127.5){0.04in}
\pscircle*(142.5,67.5){0.04in}
\pscircle(15,90){0.04in}
\pscircle(30,75){0.04in}
\pscircle(45,105){0.04in}
\pscircle(60,60){0.04in}
\pscircle(75,120){0.04in}
\pscircle*(90,45){0.04in}
\pscircle*(105,135){0.04in}
\pscircle*(120,30){0.04in}
\pscircle*(135,150){0.04in}
\pscircle(150,15){0.04in}
\pscircle(165,165){0.04in}
\rput[c](37.5,140.5){$p_{i-1}$}
\rput[c](158.5,67.5){$p_i$}
\end{pspicture}
&\rule{10pt}{0pt}&
\psset{xunit=0.01in, yunit=0.01in}
\psset{linewidth=0.005in}
\begin{pspicture}(0,0)(165,165)
\psaxes[dy=1000,dx=1000](0,0)(165,165)
\psframe[linecolor=darkgray,fillstyle=solid,fillcolor=lightgray,linewidth=0.02in](12,27)(138,160.5)
\psframe[linecolor=darkgray,fillstyle=solid,fillcolor=lightgray,linewidth=0.02in](12,57)(85.5,123)
\psline(67.5,157.5)(67.5,105)
\pscircle(82.5,97.5){0.04in}
\pscircle*(67.5,157.5){0.04in}
\pscircle*(15,90){0.04in}
\pscircle(30,75){0.04in}
\pscircle(45,105){0.04in}
\pscircle(60,60){0.04in}
\pscircle(75,120){0.04in}
\pscircle*(90,45){0.04in}
\pscircle*(105,135){0.04in}
\pscircle*(120,30){0.04in}
\pscircle*(135,150){0.04in}
\pscircle(150,15){0.04in}
\pscircle(165,165){0.04in}
\rput[c](105.5,97.5){$p_{i-1}$}
\rput[c](82.5,150.5){$p_i$}
\end{pspicture}
&\rule{10pt}{0pt}&
\psset{xunit=0.01in, yunit=0.01in}
\psset{linewidth=0.005in}
\begin{pspicture}(0,0)(195,165)
\psaxes[dy=1000,dx=1000](0,0)(195,165)
\psframe[linecolor=darkgray,fillstyle=solid,fillcolor=lightgray,linewidth=0.02in](12,49.5)(108,123)
\psframe[linecolor=darkgray,fillstyle=solid,fillcolor=lightgray,linewidth=0.02in](28,27)(168,160.5)
\psframe[linecolor=darkgray,linewidth=0.02in](12,49.5)(108,123)
\psline(30,157.5)(30,112.5)
\pscircle*(15,52.5){0.04in}
\pscircle*(30,157.5){0.04in}
\pscircle(45,90){0.04in}
\pscircle(60,75){0.04in}
\pscircle(75,105){0.04in}
\pscircle(90,60){0.04in}
\pscircle(105,120){0.04in}
\pscircle*(120,45){0.04in}
\pscircle*(135,135){0.04in}
\pscircle*(150,30){0.04in}
\pscircle*(165,150){0.04in}
\pscircle(180,15){0.04in}
\pscircle(195,165){0.04in}
\rput[c](15,37.5){$p_{i-1}$}
\rput[c](45,150.5){$p_i$}
\end{pspicture}
\end{tabular}
\end{center}
\caption{The three cases in the proof of Lemma~\ref{long-wedge}; the solid points form simple permutations.}\label{fig-long-wedge}
\end{figure}

Let us consider in detail the case of an up pin followed by a right pin.  By (W4), the left-most wedge point in $R_i\setminus R_{i-1}$ lies below $p_1$.  By separation, $p_{i-1}$ lies above $p_i$, which is itself the right-most point in $R_i$.  Therefore the wedge points in $R_i\setminus R_{i-1}$ together with $p_i$ and $p_{i-1}$ constitute a type $1$ wedge simple permutation.  The other cases follow by similar analysis; in the right-up case the wedge points in $R_i\setminus R_{i-1}$ together with $p_1$ and $p_i$ give a wedge simple permutation of type $2$, while in the left-up case a wedge simple permutation of type $2$ can be formed from the wedge points in $R_i\setminus R_{i-1}$, $p_{i-1}$, and $p_i$.
\end{proof}

We have therefore established the following theorem.

\begin{theorem}\label{sp2-really-main}
Every simple permutation of length at least $2(2048k^8)^{(2048k^8)^{2k}}$ contains a proper pin sequence of length $2k$, a parallel alternation of length $2k$, or a wedge simple permutation of length $2k$.
\end{theorem}

The proof of Theorem~\ref{sp2-main} now follows by analysing each of these cases in turn.  A parallel alternation of length $2k+2\ge 10$ contains two disjoint simple permutations of length $k$.  A type $1$ wedge simple permutation of length $2k$ contains two type $1$ wedge simple permutations of length $k$ with only one entry in common, and a type $2$ wedge simple permutation of length $2k$ contains two type $2$ wedge simple permutations of length $k$ which share two entries.  Finally, Corollary~\ref{long-pins} shows that a permutation with a proper pin sequence of length $2k+2$ contains two disjoint simple permutations of length $k$.

Brignall, Ru\v{s}kuc, and Vatter~\cite{brignall:simple-permutat:b} apply Theorem~\ref{sp2-really-main} to show that it is possible to decide whether or not a permutation class contains only finitely many simple permutations, and expatiate upon Lemma~\ref{pin-or-alt}, showing that every long simple permutation contains either a long alternation or a long ``oscillation''.

\section{Other contexts}

Although our proof is highly permutation-centric, these is no reason why analogues of Theorem~\ref{sp2-main} cannot exist for other types of object.  For example, an interval\footnote{These are also called strong intervals, partie solidaires, blocks, factors, modules, clans, congruences, and convex subsets.} in a graph is a set of vertices $X\subseteq V(G)$ such that $N(v)\setminus X=N(w)\setminus X$ for all $v,w\in X$, where $N(v)$ denotes the neighbourhood of $v$ in $G$.  A graph on $n$ vertices therefore has several trivial intervals ($\emptyset$, $V(G)$, and the singletons); a graph with no nontrivial intervals is then often called {\it prime\/} or {\it indecomposable\/} (the word simple meaning something completely different in this context).  These graphs have been the subject of considerable study, see Ehrenfeucht, Harju, and Rozenberg~\cite{ehrenfeucht:the-theory-of-2:}, Ille~\cite{ille:indecomposable-:}, and Sabidussi~\cite{sabidussi:graph-derivativ:}.

\newcommand{\A}{\mathcal{A}}
\newcommand{\B}{\mathcal{B}}
\newcommand{\dom}{\operatorname{dom}}
\newcommand{\lang}{\mathcal{L}}

The most general context for simplicity --- and thus the most general context for results such as Theorem~\ref{sp2-main} --- is relational structures.   Let $\lang$ denote a relational language (i.e., a set of relational symbols together with positive integers $n_R$ for each relational symbol $R\in\lang$, specifying the arity of $R$) and $\A$ an $\lang$-structure (i.e., a ground set $\dom(\A)$ together with interpretations of the relational symbols from $\lang$).  Following F{\"o}ldes~\cite{foldes:on-intervals-in:}, we say that the subset $X\subseteq\dom(\A)$ is an interval if the following occurs for every relation $R\in\lang$ and every $n_R$-tuple $(x_1,x_2,\dots,x_{n_R})\in\dom(\A)^{n_R}\setminus X^{n_R}$: if $x_i\in X$ then the value of $R^{\A}(x_1,x_2,\dots,x_{n_R})$ is unchanged by swapping $x_i$ with any other element of $X$.  Again the relational structure $\A$ will have the trivial intervals $\emptyset$, $\{a\}$ for all $a\in\dom(\A)$, and $\dom(\A)$ itself, and it is simple if it has no others.

In this most general context, any analogue of Theorem~\ref{sp2-main} would need to allow for more intersection between the two simple substructures.  An example demonstrating this is given in
\addtocounter{footnote}{1}%
\footnotetext{Let $\lang$ consist of a $2$-ary relation $<$ and a $k$-ary relation $R$.  Take $\A$ with $\dom(\A)=[2n]$ where $<$ is interpreted as the normal linear order on $[2n]$ and $R(1,3,5,\dots,2k-3,i)$ precisely for even $i\in[2k-2,2n]$.  This structure is simple, but all simple substructures (with at least two elements) of $\A$ must contain each of $1,3,5,\dots,2k-3$, and then to prevent these elements from containing a nontrivial interval, the simple substructure must also contain $2,4,6,\dots,2k-4$.}%
Footnote~\arabic{footnote}.

\bibliographystyle{acm}
\bibliography{../refs}

\end{document}